\documentclass[12pt]{amsart}
\usepackage{amsmath,amssymb,amsthm,amsfonts}
\usepackage{mathrsfs}
\usepackage{geometry}
\usepackage{needspace}
\usepackage[colorlinks=true,anchorcolor=blue,filecolor=blue,linkcolor=blue,urlcolor=blue,citecolor=blue]{hyperref}
\usepackage{color}
\newtheorem{theorem}{Theorem}[section]
\newtheorem{corollary}[theorem]{Corollary}
\newtheorem{lemma}[theorem]{Lemma}

\newtheorem{fact}{Fact}
\theoremstyle{definition}

\theoremstyle{remark}
\newtheorem{rmk}[theorem]{Remark}

\numberwithin{equation}{section}

\def \d{\delta}
\def \im{\operatorname{im}}
\def \K{\mathcal K}
\def \p{\partial}
\def \q{\mathfrak{q}}
\def \R{\mathbb{R}}
\def \Span{\operatorname{span}}
\def \Spec{\operatorname{Spec}}
\def \T{\wedge}

\begin{document}

\title[Spectral radius of $2$-complexes with prescribed Betti number]{Maximizing the signless Laplacian spectral radius of simplicial 2-complexes with a prescribed second Betti number}

\author[C.-M. She]{Chuan-Ming She}
\address{School of Mathematical Sciences, Anhui University, Hefei 230601, P. R. China}
\email{shecm@stu.ahu.edu.cn}

\author[Y.-Z. Fan]{Yi-Zheng Fan*}
\address{Center for Pure Mathematics, School of Mathematical Sciences, Anhui University, Hefei 230601, P. R. China}
\email{fanyz@ahu.edu.cn}
\thanks{*The corresponding author.
Supported by National Natural Science Foundation of China (Grant No. 12331012).
}

\author[Y.-M. Song]{Yi-Min Song$^\dag$}
\address{School of Mathematics and Physics, Anhui Jianzhu University, Hefei 230601, P. R. China}
\email{songym@ahjzu.edu.cn}
\thanks{$^\dag$Supported by National Natural Science Foundation of China (Grant No. 12501468).}



\begin{abstract}
We study how a prescribed second Betti number constrains the largest eigenvalue of the signless Laplacian on edges of a pure two-dimensional simplicial complex.
For each fixed positive second Betti number, we determine all maximizing complexes provided that the number of vertices is sufficiently large.
Every maximizer consists of a full cone over a complete graph together with a family of triangles avoiding the apex, any two of which share an edge; the number of added triangles equals the prescribed Betti number.
The maximizer is unique up to isomorphism except when this number is three or four, for which we describe all additional extremal complexes.
The proof uses homological constraints and Perron vector estimates to establish the full cone structure, followed by an exact Schur complement comparison to classify the added triangles.
We also obtain an asymptotic expansion of the maximum spectral radius and extend the extremal result to Betti numbers growing more slowly than the fourth root of the number of vertices.
These results provide a two-dimensional counterpart of spectral extremal theorems for connected graphs with a prescribed cyclomatic number.
\end{abstract}

\maketitle

\section{Introduction}
Spectral extremal problems ask how constraints on a combinatorial structure determine the largest eigenvalue of an associated operator.
For a connected graph $G$, its first Betti number is
$\beta_1(G)=e(G)-v(G)+1$.
Thus, prescribing the numbers of vertices and edges is equivalent to prescribing the number of vertices and the first Betti number.
This topological interpretation leads to spectral extremal problems for higher-dimensional simplicial complexes with prescribed Betti numbers.

Determining graphs with maximum adjacency spectral radius under prescribed parameters is a classical problem in spectral graph theory.
This line of research goes back to Brualdi and Hoffman \cite{brualdi1985spectral} and was developed further in \cite{brualdi1986spectral,cvetkovic1988connected}.
In particular, the maximum spectral radius of graphs with a prescribed number of edges was studied in \cite{rowlinson1988maximal,stanley1987bound}.
Related extremal questions for the signless Laplacian, also called the unoriented Laplacian, were investigated in \cite{fan2008maximizing,tam2008unoriented,zhai2020signless}.
Spectral extrema under structural restrictions have also been studied; see, for example, \cite{TAIT2017137}.
Extensions to uniform hypergraphs use adjacency tensors; Bai and Lu \cite{bai2018bound} obtained upper bounds for the spectral radius of hypergraphs with a prescribed number of edges.

An $r$-uniform hypergraph can be regarded as a pure $(r-1)$-dimensional simplicial complex by including all subsets of its edges, and conversely, the facets of a pure complex form a uniform hypergraph.
Simplicial complexes also carry homological information that is not determined by their numbers of vertices and facets alone.
Throughout this paper, homology is taken over $\R$, and $\beta_i(K)=\dim H_i(K;\R)$.
For a two-dimensional complex $K$, the second Betti number is the dimension of the space of $2$-cycles, namely $\ker\partial_2$.
It therefore depends on the linear relations among the oriented boundaries of the triangles.
An extremal problem with prescribed $\beta_2(K)$ requires control of these boundary relations as well as the distribution of faces.

We consider the signless up Laplacian on edges of a two-dimensional complex.
If $B$ is its unoriented edge--triangle incidence matrix, with $B_{e,F}=1$ when $e\subset F$ and $0$ otherwise, this operator is
$Q_1^{\rm up}(K)=BB^{\mathsf T}$.
We denote its spectral radius by $\q_1(K)$.
This is a nonnegative symmetric matrix encoding the incidence of edges and triangles, whereas the constraint on $\beta_2(K)$ involves the oriented boundary map.
The interaction between these two descriptions is central to the extremal problem considered here.
The general definitions and their relation to combinatorial Laplace operators can be found in \cite{horak2013spectra,kaufman2020high}; we recall the required notions below.

Topological constraints also arise in homeomorphic Tur\'an problems for simplicial complexes.
For a family $\mathcal F$ of $r$-dimensional complexes, let $\operatorname{ex}(n,\mathcal F)$ be the maximum number of $r$-faces in an $r$-dimensional complex on $n$ vertices that contains no member of $\mathcal F$ as a subcomplex.
When $\mathcal F$ consists of all triangulations of the sphere $S^r$, write $\operatorname{ex}(n,S^r)$.
Brown, Erd\H{o}s and S\'os \cite{SEB1973} proved that
$\operatorname{ex}(n,S^2)=\Theta(n^{5/2})$.
The same order of growth holds for every fixed closed surface, by the results of Kupavskii et al.\ \cite{KPTZ2022} in the orientable case and Sankar \cite{San2024} in the non-orientable case.
More general homeomorphic Tur\'an problems were investigated in \cite{KLNS2021,LNY2022}, while extremal problems for simplicial complexes were developed further in \cite{AGLP2025,CPS2023}.
Prescribing a Betti number fixes the dimension of a homology group, whereas excluding homeomorphs restricts the topological types of subcomplexes.
For an $r$-dimensional complex, vanishing top-dimensional homology rules out every subcomplex homeomorphic to $S^r$, but in higher dimensions the converse need not hold.

Fan, She and Zhang \cite{Fan2025hole} studied the signless Laplacian spectral radius under such topological constraints.
They proved that the full cone uniquely maximizes this spectral radius among pure $r$-dimensional complexes with vanishing top-dimensional homology, and obtained a general upper bound for complexes with a prescribed top-dimensional Betti number.
In dimension two, their bound is $\q_1(K)\le 2n-3+t$ when $\beta_2(K)=t$ and $1\le t\le n-3$, while the maximum for $t=0$ is $2n-3$.
For positive $t$, this bound does not determine the maximum or identify the extremal complexes.

In the present paper, we determine all maximizers for every fixed positive second Betti number when the number of vertices is sufficiently large.
The main structural step is to show that the homological constraint forces every maximizer to contain a full cone; a precise intersection condition then characterizes the triangles outside that cone.
For $r\ge1$, let $\K(n,r,t)$ denote the set of pure $r$-dimensional complexes on $n$ vertices with $r$-th Betti number $t$.
Let $\T_n^{r+1}$ be the full cone $r$-complex on $n$ vertices with vertex set $V$, whose facets are all $(r+1)$-subsets of $V$ containing a fixed vertex $u$.
For $r=2$, $\T_n^3=u\star K_{n-1}$, the join of $u$ and a complete graph $K_{n-1}$.
If $M\subseteq \binom{V \setminus\{u\}}{3}$, then $\T_n^3\cup M$ denotes the complex obtained by adding the triples in $M$ to the cone $\T_n^3$.

For a family $M\subseteq \binom{V\setminus\{u\}}{3}$, we say that $M$ is \emph{pairwise edge-intersecting} if any two distinct members of $M$ share a $2$-set.
We also write
\[
\T_n^{3,t}:=\T_n^3 \cup \{\{v,w,x_i\}:1\le i\le t\},
\]
where $1\le t\le n-3$ and $u,v,w,x_1,\ldots,x_t$ are distinct vertices.
Here $\{\{v,w,x_i\}:1\le i\le t\}$ is called a \emph{star pairwise edge-intersecting family}.
For an arbitrary family $M$, each added triangle, together with the three cone triangles on its edges, supports the boundary of a tetrahedron.
After choosing orientations, these boundary cycles form a basis of the second homology of $\T_n^3\cup M$, so $\beta_2(\T_n^3\cup M)=|M|$.
Our main result is the following.

\begin{theorem}\label{thm:main}
Let $t\ge 1$ be fixed. For all sufficiently large $n$, if $K\in \K(n,2,t)$ has maximum signless Laplacian spectral radius, then
\[
K=\T_n^3\cup M,
\quad |M|=t,
\]
where $M\subseteq \binom{V(K)\setminus\{u\}}{3}$ is pairwise edge-intersecting.
Conversely, every such complex is extremal and has the same signless Laplacian spectral radius as $\T_n^{3,t}$.

In particular, $\T_n^{3,t}$ is a maximizer. Moreover, the maximizer is unique up to isomorphism for $t=1,2$ and for $t\ge 5$.
For $t=3,4$, the non-star pairwise edge-intersecting subfamilies of $K_4^{(3)}$ give additional extremal complexes.
\end{theorem}

The proof also gives the asymptotic value of the maximum spectral radius.

\begin{corollary}\label{cor:asymptotic}
Let $t\ge 1$ be fixed. If $K\in \K(n,2,t)$ is extremal and $n$ is sufficiently large, then
\[
\q_1(K)=\q_1(\T_n^{3,t})=2n-3+\frac{9t}{n^3}+O_t(n^{-4}).
\]
\end{corollary}

The structural conclusion also extends to a range of growing Betti numbers.
The exact comparison within the class $\T_n^3\cup M$ works for every $1\le t\le n-3$; the restriction below comes from establishing the full cone structure.

\begin{corollary}\label{cor:growing}
Let $t=t(n)$ satisfy $3\le t\le n-3$ and $t=o(n^{1/4})$.
For all sufficiently large $n$, the complex $\T_n^{3,t}$ has maximum signless Laplacian spectral radius among all complexes in $\K(n,2,t)$.
If $t\ge 5$, it is the unique maximizer up to isomorphism.
\end{corollary}

\begin{rmk}
Let $\mathcal{G}(n,m)$ denote the set of connected graphs with $n$ vertices and $m$ edges, or equivalently, connected $1$-dimensional complexes with $n$ vertices and first Betti number $m-n+1$.
For $1\le t\le n-2$, let $\T_n^{2,t}$ be the graph obtained from the star $\T_n^2$ by adding $t$ edges among its leaves, all sharing a common vertex.
With respect to the signless Laplacian spectral radius,
$\T_n^{2,1}$ is the unique maximizing graph in $\mathcal{G}(n,n)$ for $n\ge3$ \cite{Fan},
and $\T_n^{2,2}$ is the unique maximizing graph in $\mathcal{G}(n,n+1)$ for $n\ge4$ \cite{fan2008maximizing}.
Moreover, for $n\ge5$, $\T_n^{2,3}$ is one of exactly two maximizing graphs in $\mathcal{G}(n,n+2)$ \cite{tam2008unoriented}.
For $4\le t\le n-2$, Chang and Tam \cite{chang2011} proved that $\T_n^{2,t}$ is the unique maximizing graph in $\mathcal{G}(n,n-1+t)$.
Thus, Theorem \ref{thm:main} provides a two-dimensional analogue of these results for graphs.
\end{rmk}

The proof has two stages.
First, Perron vector estimates and the constraint on the space of $2$-cycles show that all but a bounded number of triangles of an extremal complex contain a common vertex.
A face replacement argument that preserves the second Betti number then forces the full cone to be present.

Second, the triangles of the cone are indexed by the edges of $K_{n-1}$, and their adjacency graph is the triangular graph $L(K_{n-1})$.
The resolvent of its adjacency matrix has three types of entries, according as two indexing edges coincide, meet in one vertex, or are disjoint.
Taking a Schur complement reduces the remaining problem to a comparison of positive matrices whose off-diagonal entries depend only on the intersection sizes of the added triangles.
The comparison is strict unless every pair of added triangles shares an edge, yielding both the extremal value and all equality cases.

\section{Preliminaries}

\subsection{Simplicial complexes and homology}
A \emph{hypergraph} $H$ is a pair $(V,E)$, where $V=V(H)$ is the vertex set and $E=E(H)$ is the edge set, with each element of $E$ being a subset of $V$. If every edge has cardinality $r$, then $H$ is called \emph{$r$-uniform}.

An \emph{abstract simplicial complex}, or simply a \emph{complex}, $K$ on a finite set $V$ is a collection of subsets of $V$ closed under inclusion.
An \emph{$i$-face} or an \emph{$i$-simplex} of $K$ is an element of $K$ with cardinality $i+1$.
In particular, $0$-faces are called the \emph{vertices} of $K$, and $1$-faces are called the \emph{edges} of $K$.
The set of all $i$-faces is denoted by $S_i(K)$.
The dimension of a face $F$ is $|F|-1$, and the dimension of $K$ is the maximum dimension of its faces, denoted by $\dim K$.
A complex of dimension $r$ is simply called an \emph{$r$-complex}.
Faces that are maximal under inclusion are called \emph{facets}.
A complex is \emph{pure} if all its facets have the same dimension.

The \emph{join} of two complexes $K$ and $K'$ with disjoint vertex sets is denoted and defined by $K \star K'=\{F \cup F': F \in K, F' \in K'\}$.
In particular, if $K$ consists of only one vertex, say $u$, and only one nonempty face $\{u\}$,
then the join is written as $u \star K'$, which is called a \emph{cone} over $K'$.

For two $(i+1)$-faces sharing an $i$-face,  we refer to them as \emph{$i$-down neighbors}.
For two $i$-faces which are contained in a common $(i+1)$-face,
we call them \emph{ $(i+1)$-up neighbors}.
A complex $K$ is called \emph{$i$-path connected} if for any two $i$-faces $F,G$ there exists an ordering of $i$-faces $F_1,\ldots,F_m$ such that $F=F_1$, $G=F_m$, and consecutive faces are $(i+1)$-up neighbors.
For a face $F$, write $N^d(F)$ for its set of down neighbors and $N^u(F)$ for its set of up neighbors.
Given $F \in S_i(K)$ and a vertex $x \in V(K)\backslash F$, define $$N^d(F,x)=\{\{x\} \cup G \in S_i(K): G \in \p F\},$$
 where $\p F$ denotes the collection of $(i-1)$-faces of $F$.
Then we have
$$N^d(F)=\bigcup_{x \in V(K)\backslash F} N^d(F,x).$$
The \emph{degree} of an $i$-face in $K$, denoted by $d_K(F)$, is defined to be the number of $(i+1)$-faces of $K$ that contain the face $F$.

A face $F$ is \emph{oriented} if an ordering is chosen on its vertices and we write $[F]$ for the oriented face $F$.
The \emph{$i$-th chain group} $C_i(K,\R)$ is the vector space over $\R$ with a basis obtained by choosing one orientation for each $i$-face; reversing an orientation changes the sign.
We may assume that $ \emptyset \in K$ with dimension $-1$, and $C_{-1}(K,\R) \cong \R$.
For $i \ge 0$, the \emph{$i$-th boundary map} $\partial_i:C_i(K,\R)\to C_{i-1}(K,\R)$ is defined by
\[
\partial_i[v_0,\ldots,v_i]=\sum_{j=0}^i(-1)^j[v_0,\ldots,\hat{v}_{j},\ldots,v_i],
\]
 where $\hat{v}_{j}$ indicates that the vertex $ v_j $ is omitted.
These maps give rise to the chain complex:
\[
\cdots\longrightarrow C_i(K,\R)\overset{\p_i}{\longrightarrow}C_{i-1}(K,\R)
\longrightarrow\cdots\longrightarrow C_0(K,\R)
\overset{\p_0}{\longrightarrow}C_{-1}(K,\R)\longrightarrow0,
\]
which satisfies $\p_i\circ \p_{i+1}=0$ for all $i$.
The kernel $\ker\p_{i}$ is called the \emph{group of $i$-cycles}, and the image $\im\p_{i+1}$ is called the \emph{group of $i$-boundaries}.
The quotient $$H_i(K)=\ker\p_{i}/\im\p_{i+1}$$
is called the \emph{$i$-th reduced homology group of $K$}.
Its dimension, denoted by $\beta_i(K)$, is called the \emph{$i$-th Betti number} of $K$.
A nonzero element of $H_i(K)$ is called an \emph{$i$-dimensional hole} of $K$.

Let $C^i(K,\R)=\operatorname{Hom}(C_i(K,\R),\R)$ be the \emph{$i$-th cochain group} of $K$.
The \emph{$i$-th coboundary map} $\delta_i:C^i(K,\R)\to C^{i+1}(K,\R)$ is the dual of $\partial_{i+1}$, defined by $\delta_i f = f \p_{i+1}$.
Endowing $C^i(K,\mathbb{R})$ and $C^{i+1}(K,\mathbb{R})$ with positive definite inner products, respectively, we obtain
the adjoint $\delta_i^* : C^{i+1}(K,\mathbb{R}) \to C^i(K,\mathbb{R})$ of $\delta_i$, which is defined by
\[
\langle \delta_i f_1, f_2\rangle_{C^{i+1}(K,\R)} = \langle f_1, \delta_i^* f_2\rangle_{C^i(K,\R)}
\]
for all $f_1 \in C^i(K,\mathbb{R})$, $f_2 \in C^{i+1}(K,\mathbb{R})$.

Three Laplace operators \cite{horak2013spectra} on $C^i(K,\R)$ are defined as
\begin{equation}\label{Laplace}
L_i^{up}(K)=\d_{i}^*\d_{i}, ~ L_i^{down}(K)=\d_{i-1}\d_{i-1}^*, ~
L_i(K)=\d_{i}^*\d_{i} + \d_{i-1}\d_{i-1}^*
\end{equation}
which are referred to as  the \emph{$i$-th up Laplace operator}, \emph{$i$-th down Laplace operator}, and \emph{$i$-th Laplace operator} of $K$, respectively.
Eckmann \cite{eckmann1944harmonische} proved the discrete version of the Hodge theorem (see also \cite{duval2002shifted,horak2013spectra}), which can be formulated as
$\ker L_i(K) \cong H_i(K)$.

\subsection{The signless Laplacian}
The signless Laplacian of a simplicial complex was introduced in \cite{kaufman2020high} and studied further in \cite{luo2020spectrum}.
Different from the chain group generated by oriented faces, let $D_i(K,\R)$ be the vector space generated by unoriented $i$-faces, i.e. all $i$-faces of $K$.
The \emph{$i$-th signless boundary map}  $|\partial_i|:D_i(K,\R)\to D_{i-1}(K,\R)$ is defined by
\[
|\partial_i|\{v_0,\ldots,v_i\}=\sum_{j=0}^i\{v_0,\ldots,\hat{v_j},\ldots,v_i\},
\]
where $\hat{v}_j$ denotes the omission of the vertex $v_j$.
Let $D^i(K,\R)$ be the dual space of $D_i(K,\R)$ with dual basis $\{F^*: F \in S_i(K)\}$.
The \emph{$i$-th signless coboundary map} $|\d_i|: D^i(K,\R) \to D^{i+1}(K,\R)$ is defined by $ |\d_i| f =f |\p_{i+1}|$.
For each $i$, endow $D^i(K,\R)$ with a positive definite inner product that makes
 the dual basis orthonormal. (More general inner products may also be considered.)
Let $|\d_i|^*:D^{i+1}(K,\R) \to D^i(K,\R)$ be the adjoint of $|\d_i|$.
The \emph{$i$-th up signless Laplace operator} and \emph{$i$-th down signless Laplace operator} of $K$ are respectively defined by
\begin{equation}\label{eq:signless}
Q_i^{up}(K)=|\d_{i}|^* |\d_{i}|, ~ Q_i^{down}(K)=|\d_{i-1}| |\d_{i-1}|^*.
\end{equation}

Since $Q_i^{up}(K)$ has the same nonzero eigenvalues, including multiplicities, as $Q_{i+1}^{down}(K)$, we focus on discussing $Q_i^{up}(K)$ only.
For $f \in D^i(K,\R)$ and an $(i+1)$-face $\bar F$, write
\[
f(\partial \bar F):=\sum_{F\in \partial \bar F}f(F).
\]
By definition, for $f \in D^i(K,\R)$,
\begin{equation}\label{ABd}
 |\d_{i}| f = \sum_{\bar{F}\in S_{i+1}(K)} \left(\sum_{F \in \p \bar{F}} f(F) \right) \bar{F}^*=\sum_{\bar{F}\in S_{i+1}(K)} f(\partial \bar F) \bar{F}^*.
\end{equation}

Suppose that $f$ is an eigenvector of $Q_i^{up}(K)$ associated with an eigenvalue $\lambda$. Then
\begin{equation}\label{eq:edge-eig}
\lambda f(F)=(Q_i^{up}(K)f)(F)=\sum_{\bar{F}\in S_{i+1}(K):F\in\p\bar{F}}f(\p\bar{F})
=d_K(F)f(F)+\sum_{F^{\prime}\in N^{u}(F)}f(F^{\prime}).
\end{equation}
Put $g:=|\d_i|f$, so that $g(\bar F)=f(\p\bar F)$.
We have
$$ Q_{i+1}^{down}(K)g=|\d_i||\d_i|^* |\d_i|f = |\d_i|(Q_i^{up}(K)f)=\lambda |\d_i|f = \lambda g.$$
If $\lambda>0$, then $\langle g,g\rangle=\langle Q_i^{up}(K)f,f\rangle=\lambda\langle f,f\rangle>0$.
Hence $g$ is an eigenvector of $Q_{i+1}^{down}(K)$ corresponding to $\lambda$.
Thus, for $\bar{F}\in S_{i+1}(K)$,
\begin{equation}\label{eq:face-eig}
(Q_{i+1}^{down}(K)g)(\bar F)=|\bar{F}| f(\p \bar{F})+\sum_{\bar{F}' \in N^{d}(\bar{F})} f(\p \bar{F}')=\lambda f(\p \bar{F}).
\end{equation}

The matrix $Q_i^{up}(K)$ is nonnegative and symmetric. It is irreducible if and only if $K$ is $i$-path connected.
Hence, by Perron-Frobenius theory, its largest eigenvalue is the spectral radius, denoted by $\q_i(K)$.
If $Q_i^{up}(K)$ is further irreducible, then, up to a positive scalar, there exists a unique positive eigenvector associated with the spectral radius, called the \emph{Perron vector} of $Q_i^{up}(K)$.
Moreover,
\begin{equation}\label{eq:rayleigh-complex}
\q_i(K)=\max_{0\ne f\in D^i(K,\R)}\frac{\langle Q_i^{\rm up}(K)f,f\rangle}{\langle f,f\rangle}
=\max_{0\ne f\in D^i(K,\R)}\frac{\sum_{\bar F\in S_{i+1}(K)}f(\partial\bar F)^2}{\langle f,f\rangle}.
\end{equation}
In this paper, for an $r$-complex $K$, the spectral radius $\q_{r-1}(K)$ is called the \emph{signless Laplacian spectral radius} of $K$.

\section{Basic holes and extremal bounds}

Throughout this section, $r\ge1$ and $\K(n,r,t)$ denotes the family defined in the introduction.
A pure $r$-complex $K$ with $\beta_r(K)=1$ is called a \emph{basic hole} of dimension $r$ if removing any $r$-face $F$ from $K$ yields $\beta_r(K-F)=0$.
A basic hole is closely related to a hole defined as a nonzero element of the homology group.
Since $K$ is $r$-dimensional, $H_r(K)=\ker\partial_r(K)$. For a nonzero cycle $h\in H_r(K)$, let $K_h$ be the subcomplex generated by the $r$-faces with nonzero coefficients in $h$.
Clearly, $\beta_r(K_h) \ge 1$.
Moreover, $K_h$ is a basic hole if and only if $\beta_r(K_h)=1$.

Let $\Delta_{r+2}^{r+1}$ be the pure $r$-complex on $r+2$ vertices containing all $(r+1)$-subsets as facets. Let $\Diamond_{r+3}^{r+1}$ be the $r$-complex on $r+3$ vertices with facets
\[
\{\{i\}\cup F:i\in \{r+2,r+3\},\ F\in \binom{[r+1]}r\}.
\]
Both complexes are basic holes of dimension $r$ and are homeomorphic to the sphere $S^r$.

We shall use the following results from \cite[Lemma 3.1 and Corollaries 1.2--1.3]{Fan2025hole}.

\begin{theorem}\label{thm:basic-hole}
If $K$ is a basic hole of dimension $r$, then the following hold.
\begin{enumerate}
\item $K$ is $(r-1)$-path connected.
\item $d_K(F)\ge 2$ for every $(r-1)$-face $F$ of $K$.
\item $K-\bar F$ is also $(r-1)$-path connected for every $r$-face $\bar F$ of $K$.
\end{enumerate}
\end{theorem}

\begin{theorem}\label{thm:complete-skeleton-facet}
Let $K$ be a pure $r$-complex on $n$ vertices with maximum number of facets among all complexes in $\K(n,r,t)$. Then $K$ contains all possible $(r-1)$-faces and is $(r-1)$-path connected.
\end{theorem}

\begin{theorem}\label{thm:complete-skeleton-spectral}
Let $K$ be a pure $r$-complex on $n$ vertices whose signless Laplacian spectral radius $\q_{r-1}(K)$ attains the maximum among all complexes in $\K(n,r,t)$. Then $K$ contains all possible $(r-1)$-faces and is $(r-1)$-path connected.
\end{theorem}

The full cone $r$-complex $\T_n^{r+1}$ is the pure $r$-complex whose facets are all $(r+1)$-sets containing a fixed vertex. A direct calculation gives
\begin{equation}\label{cone-radius}
|S_r(\T_n^{r+1})|=\binom{n-1}{r},
\quad
\q_{r-1}(\T_n^{r+1})=rn-r^2+1.
\end{equation}
For hole-free complexes, we shall use \cite[Corollary 1.8]{Fan2025hole}.

\begin{theorem}\label{thm:holefree}
Let $K$ be a pure $r$-dimensional complex on $n$ vertices with $\beta_r(K)=0$. Then
\[
\q_{r-1}(K)\le rn-r^2+1,
\]
and equality holds if and only if $K\cong \T_n^{r+1}$.
\end{theorem}

We shall also use the rough upper bound from \cite[Theorem 1.10]{Fan2025hole}.
For a complex $K\in \K(n,r,t)$,
\begin{equation}\label{eq:rough-bound}
\q_{r-1}(K)\le rn-r^2+t+1.
\end{equation}
The cited theorem gives this bound for $1\le t\le n-r-1$.
For $t=0$, it follows from Theorem \ref{thm:holefree}; for $t\ge n-r-1$, it follows from the row-sum bound
$\q_{r-1}(K)\le(r+1)(n-r)\le rn-r^2+t+1$.

We next recall the extremal number of facets. The reduced Euler characteristic satisfies
\begin{equation}\label{eq:euler}
\widetilde\chi(K):=\sum_{i=-1}^{\dim K}(-1)^i|S_i(K)|=\sum_{i=-1}^{\dim K}(-1)^i\beta_i(K)
\end{equation}
by the Euler--Poincar\'e theorem \cite{Hatcher2002}.

\begin{lemma}\label{lem:lower-betti-zero}
Let $K$ be a pure $r$-complex on $n$ vertices with $\beta_r(K)=t$ and with maximum number of facets among all complexes in $\K(n,r,t)$. Then
\[
\beta_i(K)=0, \quad 0 \le i\le r-1.
\]
\end{lemma}

\begin{proof}
By Theorem \ref{thm:complete-skeleton-facet}, the complex $K$ contains all possible $(r-1)$-faces and is $(r-1)$-path connected.
Hence it is connected and $\beta_0(K)=0$.
Since the $(r-1)$-skeleton of $K$ is the full $(r-1)$-skeleton of an $(n-1)$-simplex, all homology groups below dimension $r-1$ vanish.
It remains only to show $\beta_{r-1}(K)=0$.

Let $\Delta_n$ be the full simplex on the same vertex set.
Take $z\in \ker \partial_{r-1}(K)$.
Since the full simplex is contractible, $z=\partial_r c$ for some $r$-chain $c$ in $\Delta_n$.
If an $r$-face $F$ appearing in $c$ is not in $K$, then adding $F$ to $K$ must create a new $r$-cycle; otherwise $K+F$ would still have $r$-th Betti number $t$ and more facets, contradicting maximality.
Thus there is an $r$-cycle $a[F]+c'$ in $K+F$, where $a\ne0$ and $c'$ is supported on $K$.
Putting $c_F=-c'/a$, we obtain $\partial_r[F]=\partial_r c_F$.
Replacing every missing face $[F]$ in $c$ by its corresponding chain $c_F$ preserves $\partial_r c=z$ and produces a chain supported on $K$.
Hence $z\in \im \partial_r(K)$, so $\beta_{r-1}(K)=0$.
\end{proof}

\begin{theorem}\label{thm:facet-max}
Let $K$ be a pure $r$-dimensional complex on $n$ vertices with $\beta_r(K)=t$. Then
\begin{equation}\label{eq:facetmax}
|S_r(K)|\le |S_r(\T_n^{r+1})|+t=\binom{n-1}{r}+t.
\end{equation}
\end{theorem}

\begin{proof}
It suffices to prove the result for a complex maximizing $|S_r(K)|$ in $\K(n,r,t)$. By Theorem \ref{thm:complete-skeleton-facet} and Lemma \ref{lem:lower-betti-zero}, $K$ contains all possible lower faces up to dimension $r-1$ and has $\beta_i(K)=0$ for $-1 \le i\le r-1$.
Using \eqref{eq:euler}, we get
\[
-1+\sum_{i=1}^r(-1)^{i+1}\binom ni+(-1)^r|S_r(K)|=(-1)^r t.
\]
 By induction on $r$ and using $\binom nr=\binom{n-1}{r}+\binom{n-1}{r-1}$, we have
\[
\sum_{i=1}^r(-1)^{i+1}\binom ni+(-1)^{r+2}\binom{n-1}{r}=1.
\]
Inequality \eqref{eq:facetmax} follows by substitution.
\end{proof}

In Theorem \ref{thm:facet-max}, the inequality \eqref{eq:facetmax} holds with equality if we take $K$ to be an $r$-complex obtained from $\T_n^{r+1}$ by adding $t$ additional $r$-faces that avoid the apex of the cone.

\section{A full cone in every extremal complex}\label{sec:structure}

In this section we prove that a spectral extremal $2$-complex with fixed positive second Betti number contains a full cone subcomplex.
After this reduction, the exact comparison in Section \ref{sec:schur} determines the extremal family.

Throughout this section, $t\ge1$ is fixed and $n$ is sufficiently large.
Let $K\in \K(n,2,t)$ be a complex maximizing $\q_1(K)$, and write $\q=\q_1(K)$.
Since $\T_n^3\cup M$ with $|M|=t$ belongs to $\K(n,2,t)$ and contains $\T_n^3$ as a proper subcomplex, by \eqref{cone-radius} and \eqref{eq:rough-bound}, we have
\begin{equation}\label{eq:qrough}
2n-3<\q\le 2n-3+t.
\end{equation}

By Theorem \ref{thm:complete-skeleton-spectral}, $K$ contains all possible edges and is $1$-path connected.
Let $f$ be the Perron vector of $Q_1^{up}(K)$, normalized so that
\begin{equation}\label{eq:normalization}
\max_{F\in S_2(K)} f(\partial F)=1.
\end{equation}
Choose $F_0=\{u,v,w\}\in S_2(K)$ with $f(\partial F_0)=1$.
By \eqref{eq:face-eig},
\[
\q=\q f(\p F_0)=3 f(\p F_0)+\sum_{F\in N^d(F_0)}f(\partial F)\le 3+|N^d(F_0)|.
\]
So, by \eqref{eq:qrough}, we have
\begin{equation}\label{low-F0nb}
|N^d(F_0)|\ge 2(n-3)+1.
\end{equation}

Let $A=V(K)\setminus F_0$, and let
\[
A_i:=\{x\in A:|N^d(F_0,x)|=i\},
\quad
A_{\le i}:=\{x\in A:|N^d(F_0,x)|\le i\}.
\]
Put $s:=|A_3|$.
Each vertex $x\in A_3$ gives a copy of $\Delta_4^{3}$ containing $F_0$.
Their boundary cycles are linearly independent, since each has a face $\{v,w,x\}$ absent from the others.
Thus $0\le s\le t$, as $\beta_2(K)=t$, and
\begin{equation}\label{A30}
|N^d(F_0)|\le 3s + 2(n-3-s) = 2(n-3)+s \le 2(n-3)+t,
\end{equation}
as each vertex of $A_3$ contributes $3$ to $|N^d(F_0)|$ and every other vertex outside $F_0$ contributes at most $2$ to $|N^d(F_0)|$.
Combining \eqref{low-F0nb} and \eqref{A30}, we have
\begin{equation}\label{eq:A3-bound}
1\le |A_3|\le t,
\quad
|A_{\le1}|\le |A_3|-1\le t-1.
\end{equation}
Indeed, if $|A_{\le1}| \ge |A_3|$, then
$$\begin{aligned}
|N^d(F_0)|&\le 3|A_3|+2|A_2|+|A_{\le 1}|\\
&=3|A_3|+2(|A|-|A_3|-|A_{\le 1}|)+|A_{\le 1}|\\
&=2|A|+|A_3|-|A_{\le 1}|\le 2|A|=2(n-3),
\end{aligned}$$
contradicting \eqref{low-F0nb}.

\begin{lemma}\label{lem:second-neighbor-lower}
With the notation above,
\[
\sum_{F_1\in N^d(F_0)}|N^d(F_1)|>4|A|^2-6t.
\]
\end{lemma}

\begin{proof}
From \eqref{eq:face-eig}, applied twice to $f(\p F_0)$, we have
\[
\q^2=\q^2 f(\p F_0) = 9+6\sum_{F_1\in N^d(F_0)}f(\partial F_1)
+\sum_{F_1\in N^d(F_0)}\sum_{F_2\in N^d(F_1)}f(\partial F_2).
\]
Since $f(\partial F)\le 1$ for every face $F$ and $|N^d(F_0)|\le 2(n-3)+t$, it follows that
\[
\begin{aligned}
\q^2 & \le  9 + 6 |N^d(F_0)| + \sum_{F_1\in N^d(F_0)} |N^d(F_1)| \\
& \le 12n-27+6t+\sum_{F_1\in N^d(F_0)}|N^d(F_1)|.
\end{aligned}
\]
Using the lower bound $\q>2n-3$ in \eqref{eq:qrough} and $|A|=n-3$ gives the result.
\end{proof}

For $x\in F_0$, define
\[
A_2^x:=\{y\in A_2:\{y\}\cup(F_0\setminus\{x\})\notin S_2(K)\}.
\]
Then $A_2=A_2^u\cup A_2^v\cup A_2^w$. We assume without loss of generality that
\[
|A_2^u|\ge |A_2^v|\ge |A_2^w|.
\]

\begin{lemma}\label{lem:aligned}
For all sufficiently large $n$,
\[
A_{\le1}=\emptyset,
\quad
A_2=A_2^u.
\]
\end{lemma}

\begin{proof}
For each $x\in A_3$, let $h_x$ be the oriented boundary cycle of the tetrahedron on $F_0\cup\{x\}$, normalized to have coefficient $1$ on the chosen oriented face $[\{v,w,x\}]$.
These cycles are linearly independent, since the face $\{v,w,x\}$ occurs in no $h_y$ with $y\ne x$.
Write
\[
Z_2(K):=\ker\partial_2(K),\qquad
W:=\Span\{h_x:x\in A_3\}.
\]
Since $K$ is two-dimensional, $\dim Z_2(K)=t$ and $\dim W=s=|A_3|$.
We call the nonzero cycles in $W$ \emph{$A_3$-holes}, and the cycles outside $W$ \emph{non-$A_3$ holes}.
Let
\[
\mathcal D:=\{F\in S_2(K):u\notin F,\ \{v,w\}\nsubseteq F\}.
\]
We claim that every cycle in $Z_2(K)$ whose coefficients on $\mathcal D$ vanish belongs to $W$.
Indeed, its only possible faces avoiding $u$ are $\{v,w,x\}$.
If such a face has nonzero coefficient, cancellation along the edges $\{v,x\}$ and $\{w,x\}$ forces both $\{u,v,x\}$ and $\{u,w,x\}$ to be present, so $x\in A_3$.
Subtracting the appropriate linear combination of the cycles $h_x$ leaves a cycle supported entirely on faces containing $u$.
Such a cycle is zero: each edge avoiding $u$ belongs to at most one face in its support.
This proves the claim.

Starting with $K_0=K$, suppose that $Z_2(K_j)\ne W$.
Choose a cycle in $Z_2(K_j)\setminus W$ and a face $F_j\in\mathcal D$ on which its coefficient is nonzero, and put $K_{j+1}=K_j-F_j$, retaining the proper faces of $F_j$.
Then
\[
Z_2(K_{j+1})=\{z\in Z_2(K_j):z_{F_j}=0\}
\]
has codimension one in $Z_2(K_j)$ and still contains $W$.
Moreover, each edge of $F_j$ lies in another face of the chosen cycle, so this deletion preserves purity.
After exactly $t-s$ steps we obtain a complex $K'$ with $Z_2(K')=W$.
None of the deleted faces lies in $\{F_0\}\cup N^d(F_0)$, so $N^d(F_0)$ and all the sets $A_i$ are unchanged.
Each deleted face contributes at most four to $\sum_{F_1\in N^d(F_0)}|N^d(F_1)|$:
a face contained in $A$ contributes zero, while a face $\{v,x,y\}$ can be adjacent only to $\{u,v,x\}$, $\{u,v,y\}$, $\{v,w,x\}$ and $\{v,w,y\}$ in $N^d(F_0)$, and the case $\{w,x,y\}$ is similar.
Thus the total contribution of the deleted faces is at most $4(t-s)$.

We record a local consequence used in the counting below.
For distinct $x,y\in A_2\cup A_3$, put $s_{xy}=|\{x,y\}\cap A_3|$, and let $k$ be the number of faces among $\{u,x,y\}$, $\{v,x,y\}$ and $\{w,x,y\}$ present in $K'$.
The subcomplex induced on $F_0\cup\{x,y\}$ has $5+s_{xy}+k$ triangles.
Its cycle space is contained in $\Span\{h_z:z\in\{x,y\}\cap A_3\}$ and thus has dimension at most $s_{xy}$.
The boundary matrix of any complex on five vertices has rank at most $\binom42=6$, since its image lies in the six-dimensional cycle space of $K_5$.
Consequently $5+s_{xy}+k-6\le s_{xy}$, so
\begin{equation}\label{eq:local-extra-face}
\bigl|\{\{u,x,y\},\{v,x,y\},\{w,x,y\}\}\cap S_2(K')\bigr|\le1.
\end{equation}
Moreover, no tetrahedron boundary in $K'$ can contain a vertex of $A_{\le1}$, since every cycle in $W$ is supported on $F_0\cup A_3$.
Hence $|N^d(F_1,y)|\le2$ whenever $y\notin F_1$ and $(F_1\cup\{y\})\cap A_{\le1}\ne\emptyset$.

A straightforward partition of the set $N^d(F_0)$ gives
\begin{align}\label{eq:splitS}
\sum_{F_1\in N^d(F_0)}|N^d(F_1)|
& =
\sum_{x\in A_3}\sum_{F_1\in N^d(F_0,x)}|N^d(F_1)|
+\sum_{x\in A_{\le1}}\sum_{F_1\in N^d(F_0,x)}|N^d(F_1)|\notag\\
&\quad+
\sum_{s\in F_0}\sum_{x\in A_2^s}\sum_{F_1\in N^d(F_0,x)}|N^d(F_1)|.
\end{align}
We now discuss the upper bounds for the three summands in \eqref{eq:splitS} in the complex $K'$.

\begin{fact}\label{eq:S-A3}
In the complex $K'$,
$$\sum_{x\in A_3}\sum_{F_1\in N^d(F_0,x)}|N^d(F_1)|\le |A_3|(4|A|+|A_{3}|+2|A_{\le 1}|+4).$$
\end{fact}

\begin{proof}
Note that for each $x \in A_3$, $N^d(F_0,x)=\{\{u,v,x\},\{u,w,x\},\{v,w,x\}\}$.
We express the sum $S:=\sum_{x\in A_3}\sum_{F_1\in N^d(F_0,x)}\left|N^d(F_1)\right|$ as
\begin{equation}\label{NdF0A31}
\begin{aligned}
S &=\sum_{x\in A_3}\sum_{F_1\in N^d(F_0,x)}\sum_{y \in F_0\setminus F_1}|N^d(F_1,y)|+\sum_{x\in A_3}\sum_{F_1\in N^d(F_0,x)}\sum_{y \notin F_0\cup\{x\}}|N^d(F_1,y)|\\
&=\sum_{x\in A_3}\left(|N^d(\{u,v,x\},w)|+|N^d(\{u,w,x\},v)|
   +|N^d(\{v,w,x\},u)|\right)\\
& +\sum_{x\in  A_3}\sum_{F_1\in N^d(F_0,x)}\sum_{y\in  A_2}|N^d(F_1,y)|
+\sum_{x\in  A_3}\sum_{F_1\in N^d(F_0,x)}\sum_{y\in ( A_3\cup  A_{\le 1})\backslash \{x\}}|N^d(F_1,y)|.
\end{aligned}
\end{equation}
	
For any $x\in  A_3$, we have \[|N^d(\{u,v,x\},w)|=|N^d(\{u,w,x\},v)|=|N^d(\{v,w,x\},u)|=3.\]
Meanwhile, for every $x\in  A_3$,
\[\sum_{F_1\in N^d(F_0,x)}\sum_{y\in  A_2}|N^d(F_1,y)|=\sum_{y\in  A_2}\sum_{F_1\in N^d(F_0,x)}|N^d(F_1,y)|\le 4| A_2|,
\]
where the last inequality holds because $\sum_{F_1\in N^d(F_0,x)}|N^d(F_1,y)|\le 4$ for any $y\in A_2$ by the following reason.
Suppose $y \in A_2^u$. By \eqref{eq:local-extra-face}, $K'$ contains at most one of $\{y,u,x\}$, $\{y,v,x\}$ and $\{y,w,x\}$. The faces containing two vertices of $F_0$ contribute two to the sum, and this additional face contributes at most two, giving the bound four.

Similarly, for each $x\in  A_3$,
\[\begin{aligned}
&\sum_{F_1\in N^d(F_0,x)}\sum_{y\in ( A_3\cup  A_{\le 1})\backslash \{x\}}|N^d(F_1,y)|\\
&=\sum_{y\in ( A_3\cup  A_{\le 1})\backslash \left\{x\right\}}\sum_{F_1\in N^d(F_0,x)}|N^d(F_1,y)|\\
& = \sum_{y\in  A_3\backslash \{x\}} \sum_{F_1\in N^d(F_0,x)}|N^d(F_1,y)| +
\sum_{y\in  A_{\le 1}}  \sum_{F_1\in N^d(F_0,x)}|N^d(F_1,y)|\\
& \le 5 (| A_3|-1) + 6 | A_{\le 1}|.
\end{aligned}\]
Here the bound five for $y\in A_3\setminus\{x\}$ follows from \eqref{eq:local-extra-face}, and the bound six for $y\in A_{\le1}$ follows from the preceding observation.

Therefore, by \eqref{NdF0A31}, we obtain
$$
\begin{aligned}
\sum_{x\in  A_3}\sum_{F_1\in N^d(F_0,x)}\left|N^d(F_1)\right|
&\le 9| A_3|
  +\sum_{x\in A_3}\left(4| A_2|+5(| A_3|-1)+6| A_{\le 1}|\right)\\
&=| A_3|(4| A|+| A_{3}|+2| A_{\le 1}|+4).
\end{aligned}
$$
\end{proof}

\begin{fact}\label{eq:S-Ale1}
In the complex $K'$,
$$\sum_{x\in A_{\le 1}}\sum_{F_1\in N^d(F_0,x)}|N^d(F_1)|\le
 2|A||A_{\le 1}|.$$
\end{fact}

\begin{proof}
For each $x\in A_{\le 1}$, $|N^d(F_0,x)|\le 1$.
If $N^d(F_0,x) =:\{F_1\}$, then for each vertex $y \notin F_1$, $|N^d(F_1,y)| \le 2$; otherwise we would have a non-$A_3$-hole $\Delta_4^3$ of $K'$, a contradiction.
So, $$|N^d(F_1)|=\sum_{y \notin F_1} |N^d(F_1,y)| \le 2(n-3)=2|A|.$$
The inequality follows.
\end{proof}

\begin{fact}\label{eq:S-A2}
In the complex $K'$,
$$\sum_{x\in F_0}\sum_{y\in A_2^x}\sum_{F_1\in N^d(F_0,y)}|N^d(F_1)|\le 4|A||A_2|.$$
Moreover, if $A_2^v\ne \emptyset$, then
$$\sum_{x\in F_0}\sum_{y\in A_2^x}\sum_{F_1\in N^d(F_0,y)}|N^d(F_1)|\le 4|A||A_2|-2(|A_2|-1).$$
\end{fact}

\begin{proof}
By symmetry, it suffices to analyze the case of $x=u$.
Given $y \in A_2^u$, $N^d(F_0,y)=\{\{y,u,v\},\{y,u,w\}\}$, and
\begin{equation}\label{NdF0A21}
\begin{aligned}
\sum_{F_1\in  N^d(F_0,y)}|N^d(F_1)|&=|N^d(\{y,u,v\})|+|N^d(\{y,u,w\})|\\
&=|N^d(\{y,u,v\},w)|+|N^d(\{y,u,w\},v)|\\
&+ \sum_{z\in A_3}(|N^d(\{y,u,v\},z)|+|N^d(\{y,u,w\},z)|)\\
& +\sum_{z\in A_{\le 2}\backslash \{y\}}\sum_{F_1\in N^d(F_0,y)}|N^d(F_1,z)|\\
& \le 4+4|A_3| +\sum_{z\in A_{\le 2}\backslash \{y\}}\sum_{F_1\in N^d(F_0,y)}|N^d(F_1,z)|.
\end{aligned}
\end{equation}

We bound the term $\sum_{z\in A_{\le 2}\backslash \{y\}}\sum_{F_1\in N^d(F_0,y)}|N^d(F_1,z)|$ in \eqref{NdF0A21} according to the partition
$A_{\le 2}\backslash \{y\}=(A_2^u\backslash\{y\}) \cup (A_2^v\cup A_2^w) \cup A_{\le 1}$.

\emph{Case} 1:  $z\in A_2^u \backslash \{y\}$.
We have
 $$\sum_{F_1\in N^d(F_0,y)}|N^d(F_1,z)|=|N^d(\{y,u,v\},z)|+|N^d(\{y,u,w\},z)|\le 4.$$

\emph{Case} 2:  $z\in A_2^v\cup A_2^w$.
We have
  $\sum_{F_1\in N^d(F_0,y)}|N^d(F_1,z)|\le 3.$
This follows from \eqref{eq:local-extra-face}: the faces containing two vertices of $F_0$ contribute one to the sum, and at most one additional face contributes at most two.

\emph{Case} 3: $z\in A_{\le 1}$. In this case, we have
 $\sum_{F_1\in N^d(F_0,y)}|N^d(F_1,z)|\le 4.$

Combining the above three cases, we have
\begin{equation}\label{NdF0A22}
\begin{aligned}
\sum_{z\in A_{\le 2}\backslash \{y\}}\sum_{F_1\in N^d(F_0,y)}|N^d(F_1,z)|
& \le 4(|A_2^u|-1)+3(|A_2^v|+|A_2^w|)+4|A_{\le 1}|\\
&=4|A_{\le 2}|-(|A_2^v|+|A_2^w|)-4.
\end{aligned}
\end{equation}

By \eqref{NdF0A21} and \eqref{NdF0A22}, we have
\begin{equation}\label{NdF0A23}
\begin{aligned}
\sum_{y\in A_2^u}\sum_{F_1\in N^d(F_0,y)}|N^d(F_1)|
& \le \sum_{y\in A_2^u}\left(4+4|A_3|+4|A_{\le 2}|-(|A_2^v|+|A_2^w|)-4\right)\\
&=\sum_{y\in A_2^u}\left(4|A|-(|A_2^v|+|A_2^w|)\right).\\
\end{aligned}
\end{equation}
	
By a similar discussion, we can derive that
\begin{equation}\label{NdF0A24}
\begin{aligned}
\sum_{y\in A_2^v}\sum_{F_1\in N^d(F_0,y)}|N^d(F_1)| & \le \sum_{y\in A_2^v}\left(4|A|-(|A_2^u|+|A_2^w|)\right),\\
\sum_{y\in A_2^w}\sum_{F_1\in N^d(F_0,y)}|N^d(F_1)|& \le \sum_{y\in A_2^w}\left(4|A|-(|A_2^u|+|A_2^v|)\right).
\end{aligned}
\end{equation}
Combining \eqref{NdF0A23} and \eqref{NdF0A24}, we get
\begin{equation}\label{NdF0A25}
\begin{aligned}
\sum_{x\in F_0}\sum_{ y \in A_2^x}\sum_{F_1\in N^d(F_0,y)}|N^d(F_1)|
 &\le 4|A||A_2|-2(|A_2^u||A_2^v|
 +|A_2^u||A_2^w|+|A_2^v||A_2^w|)\\
&=4|A||A_2|-|A_2|^2+|A_2^u|^2+|A_2^v|^2+|A_2^w|^2\\
& \le 4|A||A_2|.
\end{aligned}
\end{equation}

If $A_2^v\ne \emptyset$, then
\[|A_2^u|^2+|A_2^v|^2+|A_2^w|^2\le |A_2^u|^2+(|A_2^v|+|A_2^w|)^2 \le(|A_2|-1)^2+1.\]
Substituting back into \eqref{NdF0A25}, we get the desired upper bound.
\end{proof}

Now we return to the upper bound for $\sum_{F_1\in N^d(F_0)}|N^d(F_1)|$.
Note that the above Facts are discussed for the complex $K'$, which is obtained from $K$ by deleting some $2$-faces, and those deleted faces contribute at most $4(t-|A_3|)$ to $\sum_{F_1\in N^d(F_0)}|N^d(F_1)|$.
So, by Facts \ref{eq:S-A3}, \ref{eq:S-Ale1} and \ref{eq:S-A2},
\begin{equation}\label{eq:S-total-F}
\begin{aligned}
\sum_{F_1\in N^d(F_0)}|N^d(F_1)|&\le|A_3|(4|A|+|A_{3}|+2|A_{\le 1}|+4)+2|A||A_{\le 1}|\\
&\quad+4|A||A_2|+4(t-|A_3|)\\
&=4|A|^2-2|A||A_{\le 1}|+|A_3|(|A_{3}|+2|A_{\le 1}|)+4t.
\end{aligned}
\end{equation}
If $A_2^v \ne \emptyset$, then
\begin{equation}\label{eq:S-A2-improved}
\sum_{F_1\in N^d(F_0)}|N^d(F_1)| \le 4|A|^2-2|A||A_{\le 1}|+|A_3|(|A_{3}|+2|A_{\le 1}|)+4t - 2(|A_2|-1).
\end{equation}

By Lemma \ref{lem:second-neighbor-lower} and \eqref{eq:S-total-F}, we have
\[
2|A||A_{\le1}|\le |A_3|(|A_3|+2|A_{\le1}|)+10t.
\]
Since $|A_3| \le t$ and $|A_{\le1}|\le t-1$ by \eqref{eq:A3-bound}, we have
$$ 2|A||A_{\le 1}|\le 3t^2+8t.$$
If $A_{\le 1}\ne \emptyset$, then
 $$2(n-3)=2|A|\le 2|A||A_{\le 1}|\le 3t^2+8t.$$
This yields a contradiction for sufficiently large $n$.
Thus, $A_{\le 1}=\emptyset$.

Now assume $A_2^v\ne\emptyset$.
By Lemma \ref{lem:second-neighbor-lower} and \eqref{eq:S-A2-improved}, noting that $A_{\le 1}=\emptyset$, we have
$$ 0= 2|A||A_{\le 1}|\le 3t^2+8t - 2(|A_2|-1).$$
As $|A_2|=|A|-|A_3| \ge n-3-t$,
we also get a contradiction for sufficiently large $n$.
Therefore, $A_2^v=A_2^w=\emptyset$, and hence $A_2=A_2^u$.
\end{proof}

Let
\[
M_u:=\{F\in S_2(K):u\notin F\},
\quad
M_u(vw):=\{F\in M_u:\{v,w\}\subseteq F\}.
\]

\begin{lemma}\label{lem:bounded-Mu}
For all sufficiently large $n$,
\[
|M_u|<5t^2+10t.
\]
\end{lemma}

\begin{proof}
Let
\[
S:=\sum_{F_1\in N^d(F_0)}|N^d(F_1)|
=|\{(F_1,F_2):F_1\in N^d(F_0),\ F_2\in N^d(F_1)\}|.
\]
By Lemma \ref{lem:aligned}, $A_{\le 1}=\emptyset$ and $V(K) \setminus F_0=A=A_2^u \cup A_3$.
It is easy to see that  $|M_u(vw)|=|A_3|\le t$, and each face $F \in M_u(vw)$ has the form $F=\{v,w,x\}$ for some vertex $x \in A_3$.

Let $K'$ be the complex defined in Lemma \ref{lem:aligned}, which is obtained from $K$ by
deleting $t-|A_3|$ $2$-faces that avoid $u$ and contain at most one of $v,w$, whose contribution to $S$ is at most $4(t-|A_3|)$.
By double counting, we count pairs $(F_1,F_2)$ in $K'$ such that $F_2 \in S_2(K')$ and $F_1 \in N^d(F_0) \cap N^d(F_2) \cap S_2(K')$.

\emph{Case} 1: $u \notin F_2$ (or equivalently, $F_2\in M_u$).

\begin{enumerate}

\item[(1.1)]  If $F_2\in M_u(vw)$, then $F_2=\{v,w,x\}$ for some vertex $x \in A_3$, and it contributes $|A_3|+1$ choices for $F_1$.

\item[(1.2)] If $|F_2\cap\{v,w\}|=1$ and $|F_2\cap A_3|=2$, then $F_2$ contributes  $4$ choices for $F_1$.

\item[(1.3)] If $|F_2\cap\{v,w\}|=1$ and $|F_2\cap A_3| \le 1$, then $F_2$ contributes at most $3$ choices for $F_1$.

\item[(1.4)] If $F_2 \cap \{v,w\} = \emptyset$, then $F_2 \subseteq A$ and no choice exists for $F_1$.

\end{enumerate}
		
\Needspace{5\baselineskip}
\emph{Case} 2: $u\in F_2$.

\begin{enumerate}

\item[(2.1)] If $F_2=F_0=\{u,v,w\}$, then $F_2$ contributes $2|A|+|A_3|$ choices for $F_1$.

\item[(2.2)] If $F_2 \in \{\{x,u,v\},\{x,u,w\}\}$ for some $x\in A_3$, then $F_2$  contributes $|A|+1$ choices for $F_1$.

\item[(2.3)] If $F_2\in \{\{x,u,v\},\{x,u,w\}\}$ for some $x\in A_{2}$, then $F_2$  contributes $|A|$ choices for $F_1$.

\item[(2.4)] If $F_2=\{x,y,u\}$ for two distinct vertices $x,y\in A$, then $F_2$  contributes $4$ choices for $F_1$.

\end{enumerate}

Note that $K$ has at most $\binom{n-1}{2}+t$ facets by Theorem \ref{thm:facet-max}, among which $|M_u|$ facets avoid $u$, and $1+2|A|$ facets contain the edge $\{u,v\}$ or $\{u,w\}$.
So, the number of faces $F_2$ in Case 2.4 is bounded above by
$$\binom{n-1}{2}+t-|M_u|-2|A|-1.$$
Now, returning to $K$, we have
\[\begin{aligned}
\sum_{F_1\in N^d(F_0)}|N^d(F_1)|& \le |A_3|(|A_3|+1)+2\cdot \binom{|A_3|}{2}\cdot 4+3|M_u|\\
& ~~~ +(2|A|+|A_3|)+2|A_3|(|A|+1) +2|A_2||A|\\
&~~~ +4\left(\binom{n-1}{2}+t-|M_u|-2|A|-1\right)+4(t-|A_3|)\\
&=4|A|^2+|A_3|^2+8\binom{|A_3|}{2}+8t-|M_u|,
\end{aligned}\]
Combining this with Lemma \ref{lem:second-neighbor-lower} and $|A_3|\le t$, we have
\[
|M_u|<|A_3|^2+8\binom{|A_3|}{2}+14t\le 5t^2+10t.
\]
\end{proof}

\begin{lemma}\label{lem:cone-completion}
For all sufficiently large $n$, the extremal complex $K$ contains every $2$-face containing $u$. Hence
\[
K=\T_n^3\cup M_u.
\]
\end{lemma}

\begin{proof}
Put $M=5t^2+10t$. By Lemma \ref{lem:bounded-Mu}, $|M_u|<M$. If $F\in M_u$, then by \eqref{eq:face-eig},
\[
(\q-3)f(\partial F)=\sum_{F'\in N^d(F)}f(\partial F')\le |M_u|-1+3<M+2.
\]
Using $\q>2n-3$ from \eqref{eq:qrough}, we obtain
\begin{equation}\label{eq:small-noncone}
f(\partial F)<\frac{M+2}{2(n-3)}, \text{~ for~} F\in M_u.
\end{equation}
Next, by \eqref{eq:face-eig} applied to $F_0$,
\[
\q-3=\sum_{F\in N^d(F_0)}f(\partial F)
=\sum_{F\in M_u(vw)}f(\partial F)+\sum_{F\in N^d(F_0)\setminus M_u(vw)}f(\partial F).
\]
Together with \eqref{eq:small-noncone}, this gives
\[
\sum_{F\in N^d(F_0)\setminus M_u(vw)}f(\partial F)
>2(n-3)-\frac{M(M+2)}{2(n-3)}.
\]
Since $N^d(F_0)\setminus M_u(vw)$ has $2(n-3)$ faces, and for each face $F \in N^d(F_0)\setminus M_u(vw)$, $f(\p F) \le 1$, we have
\begin{equation}\label{eq:large-cone-boundary}
f(\partial F)>1-\frac{M(M+2)}{2(n-3)}.
\end{equation}

For any edge $G$ not containing $u$, the equation \eqref{eq:edge-eig} gives
\begin{equation}\label{eq:small-edge}
f(G)=\frac{1}{\q}\sum_{F\supseteq G}f(\partial F)
<\frac{M+1}{2n-3}.
\end{equation}
Similarly,
\[
f(\{u,v\})\le \frac{n-2}{\q}<\frac{n-2}{2n-3}.
\]
For every $x\in A$, the face $\{u,v,x\}$ lies in $N^d(F_0)\setminus M_u(vw)$; hence by \eqref{eq:large-cone-boundary} and \eqref{eq:small-edge},
\begin{equation}\label{eq:ux-large}
f(\{u,x\})=f(\partial\{u,v,x\})-f(\{u,v\})-f(\{v,x\})
>\frac12-\frac{(M+2)^2}{2(n-3)}.
\end{equation}

Suppose, for a contradiction, that $C=\{u,x,y\}\notin S_2(K)$ for some $x,y \in A$.
Consider $K^+=K+C$. If $\beta_2(K^+)=t$, then by \eqref{eq:rayleigh-complex}, the Rayleigh quotient strictly increases after adding $C$, since $f(\partial C)>0$, contradicting the extremality of $K$. Hence $\beta_2(K^+)=t+1$, and $K^+$ contains a new basic hole $H$ containing $C$: choose a cycle using $C$ with inclusion-minimal support.
The hole $H$ must contain a face $F \in M_u$ that avoids $u$; otherwise, $H$ would be a cone with $\beta_2(H)=0$, a contradiction.
Now, replacing $F$ by $C$, we get a complex $\widetilde K:=K^+-F=K-F+C$.
A nonzero cycle supported on the basic hole $H$ has nonzero coefficient on $F$.
Thus the coordinate condition $z_F=0$ has codimension one in $Z_2(K^+)$, and $\beta_2(\widetilde K)=\beta_2(K^+)-1=t$.
By Theorem \ref{thm:basic-hole}, each edge of $F$ lies in another face of $H$.
Hence $\widetilde K$ is pure and has the same edges as $K$, so $\widetilde K\in\K(n,2,t)$ and the same vector $f$ can be used in the Rayleigh quotient.
Using the Rayleigh quotient \eqref{eq:rayleigh-complex}, \eqref{eq:small-noncone}, and \eqref{eq:ux-large}, for sufficiently large $n$,
\begin{align*}
(\q_1(\widetilde K)-\q)\langle f,f\rangle
&\ge f(\partial C)^2-f(\partial F)^2\\
&\ge \left(2\left(\frac12-\frac{(M+2)^2}{2(n-3)}\right)\right)^2
-\left(\frac{M+2}{2(n-3)}\right)^2>0.
\end{align*}
This contradicts the extremality of $K$. Hence $K$ contains the cone $\T_n^3$.
\end{proof}

\begin{theorem}\label{thm:structural}
Let $t\ge1$ be fixed. For all sufficiently large $n$, every complex $K\in \K(n,2,t)$ with maximum signless Laplacian spectral radius has the form
\[
K=\T_n^3\cup M,
\quad
M\subseteq \binom{V(K)\setminus\{u\}}{3},
\quad
|M|=t.
\]
Moreover,
\[
\q_1(K)=2n-3+\frac{9t}{n^3}+O_t(n^{-4}).
\]
\end{theorem}

\begin{proof}
The structural statement follows from Lemma \ref{lem:cone-completion}. Since $\T_n^3$ has $\binom{n-1}{2}$ facets, Theorem \ref{thm:facet-max} gives $|M|\le t$. Conversely, the cone has second Betti number zero and adding one face increases this number by at most one, so $t=\beta_2(K)\le|M|$. Hence $|M|=t$.

It remains to compute the asymptotic value. Let $C$ be the set of cone faces and let $M$ be the set of non-cone faces.
For every cone face $F\in C$, by $f(\p F_0)=\max f(\p F)=1$, \eqref{eq:large-cone-boundary}, and \eqref{eq:ux-large}, one has
$$f(\partial F)=1-O_t(n^{-1}),$$
whence
\[
\sum_{F\in C}f(\partial F)=\binom{n-1}{2}-O_t(n).
\]
For every $F\in M$, \eqref{eq:small-noncone} gives
\[
f(\partial F)=O_t(n^{-1}).
\]
For every edge $G\subseteq V\setminus\{u\}$, by \eqref{eq:edge-eig},
\[
f(G)=\frac{1}{\q}\biggl(f(\partial(\{u\}\cup G))+\sum_{F\in M:G\subseteq F}f(\partial F)\biggr)
=\frac1{2n}+O_t(n^{-2}).
\]
Thus, for every $F\in M$,
\begin{equation}\label{f-M-value}
f(\partial F)=\frac{3}{2n}+O_t(n^{-2}).
\end{equation}

Summing the $Q_2^{\rm down}$ eigenvalue equation \eqref{eq:face-eig} over all cone faces gives the exact identity
\begin{equation}\label{eq:cone-sum-identity}
(\q-(2n-3))\sum_{F\in C}f(\partial F)=3\sum_{F\in M}f(\partial F).
\end{equation}
Indeed, each cone face has $2(n-3)$ cone down-neighbors, and each non-cone face has exactly three cone down-neighbors. Therefore
\[
\q-(2n-3)=\frac{3t(\frac{3}{2n}+O_t(n^{-2}))}{\binom{n-1}{2}-O_t(n)}
=\frac{9t}{n^3}+O_t(n^{-4}).
\]
\end{proof}

\begin{rmk}
The structural proof above is intentionally separated into two parts. Lemmas \ref{lem:aligned} and \ref{lem:bounded-Mu} show that all but $O_t(1)$ faces contain a common vertex. Lemma \ref{lem:cone-completion} is a short Perron-vector switching argument showing that once the exceptional part is bounded, all cone faces must be present.
\end{rmk}

\section{Spectral comparison within the cone class}\label{sec:schur}

In this section we determine which $3$-uniform family $M$ maximizes $\q_1(\T_n^3\cup M)$ after the structural reduction. This part is exact and does not require $t$ to be fixed.

Since $Q_1^{up}$ and $Q_2^{down}$ have the same nonzero eigenvalues, we may work with
$Q_2^{down}$.
By definition, for a pure $2$-complex $K$,
$$ Q_2^{down}(K) = 3 I + A_2^{down}(K),$$
where $A_2^{down}(K)$ is the down-adjacency matrix of facets of $K$, defined by
\[
A_2^{down}(K)_{F,G}=\begin{cases}
1,&F\ne G\text{ and }F,G\text{ share an edge},\\
0,&\text{otherwise}.
\end{cases}
\]

Let $N=n-1$ and let $u$ be the apex of $\T_n^3$.
The cone faces of $\T_n^3$ are naturally indexed by the edges of the complete graph $K_N$ on $V\setminus\{u\}$.
Let $L(K_N)$ be the line graph of $K_N$, and let $A(L(K_N))$ be the adjacency matrix of $L(K_N)$.
For a triple $F\subseteq V\setminus\{u\}$, let $b_F$ be the incidence (column) vector of the three boundary edges of $F$ in $E(K_N)$, that is, $(b_F)_e=1$ if $e \in \p F$ and $(b_F)_e=0$ otherwise.

Now taking $K=\T_n^3\cup M$, we have the following block form for $A_2^{down}(K)$:
\begin{equation}\label{eq:blockAM}
A_M: = A_2^{down}(\T_n^3\cup M)=
\begin{pmatrix}
A_0&B\\
B^T&C_M
\end{pmatrix}
\end{equation}
where $A_0=A(L(K_N))$ is the principal submatrix of $A_2^{down}(\T_n^3\cup M)$ on the cone faces indexed by the edges of $K_N$, the columns of $B$ are the vectors $b_F$ for $F\in M$, and $C_M$ is the down-adjacency matrix of faces of $M$.

For $N\ge4$, the spectrum of $A(L(K_N))$ is
\begin{equation}\label{eq:triangular-spectrum}
\Spec A(L(K_N))=\{(2N-4)^{(1)}, (N-4)^{(N-1)}, (-2)^{\left(\binom{N}{2}-N\right)}\},
\end{equation}
where $\lambda^{(m)}$ means $\lambda$ as an eigenvalue with multiplicity $m$.
Let $\rho_0:=2N-4=2n-6$, and let $\rho(A)$ denote the spectral radius (or the largest eigenvalue) of a nonnegative symmetric matrix $A$.
Then
\[
\q_1(\T_n^3\cup M)=3+\rho(A_M).
\]

\begin{lemma}\label{lem:exact-comparison}
Fix $n$ and $1 \le t\le n-3$. Among all complexes of the form
\[
\T_n^3\cup M,
\quad
M\subseteq \binom{V\setminus\{u\}}{3},
\quad
|M|=t,
\]
the signless Laplacian spectral radius is maximized when the triples in $M$ are pairwise edge-intersecting. In particular, the star family
\[
M^*=\{\{v,w,x_i\}:1\le i\le t\}
\]
is a maximizer.
\end{lemma}

\begin{proof}
If $n=4$, then $t=1$ and there is only one possible family $M$, so the conclusion is immediate.
We may therefore assume $n\ge5$.
Let $\lambda_M=\rho(A_M)$.
Since $\T_n^3\cup M$ is $1$-path connected, $A_M$ is irreducible and has $A_0$ as a proper principal submatrix. Hence $\lambda_M>\rho_0$.
Fix $\lambda>\rho_0$ and put
\[
R_\lambda=(\lambda I-A_0)^{-1}.
\]
Since $\lambda I-A_0$ is an invertible and irreducible $M$-matrix, $R_\lambda$ is entrywise positive; see \cite[Chapter 6, Theorem 2.7]{BerP}.

We next derive \eqref{eq:schur-root}.
Suppose first that $\lambda>\rho_0$ is an eigenvalue of $A_M$ with eigenvector written according to the block decomposition as $(x,y)^\top$, where $x$ is indexed by cone faces and $y$ is indexed by the faces of $M$.
The eigenvalue equation is
\[
A_0x+By=\lambda x,
\quad
B^Tx+C_My=\lambda y.
\]
The first equation gives $x=R_\lambda B y$, implying $y \ne 0$.
Substituting this into the second equation yields
\[
(C_M+B^TR_\lambda B)y=\lambda y.
\]
Equivalently,
\[
(\lambda I-C_M-B^\top R_\lambda B)y=0.
\]
Thus, $\lambda$ is an eigenvalue of $A_M$ above $\rho_0$ if and only if the Schur complement
\[
\lambda I-C_M-B^TR_\lambda B
\]
is singular.

When $\lambda=\lambda_M$, choose $(x,y)^\top$ to be the Perron vector of $A_M$. Then $y$ is positive.
Moreover, since $R_\lambda$ is positive and each column $b_F$ of $B$ has three positive entries,  the matrix $B^TR_{\lambda_M}B$ is positive.
Hence $C_M+B^TR_{\lambda_M}B$ is a positive matrix.
By the Perron-Frobenius theorem, $\lambda_M$ is exactly the spectral radius of the smaller matrix $C_M+B^TR_{\lambda_M}B$.
Therefore $\lambda_M$ is characterized by the equation
\begin{equation}\label{eq:schur-root}
\lambda=\rho(C_M+B^T R_\lambda B),
\end{equation}
with $\lambda>\rho_0$.
Indeed, this equation has a unique solution above $\rho_0$:
since $R_\lambda'=-R_\lambda^2$ and $R_\lambda$ is entrywise positive, every entry of $B^TR_\lambda B$ strictly decreases with $\lambda$. Perron-Frobenius theory implies that $\rho(C_M+B^T R_\lambda B)$ strictly decreases, whereas the left side $\lambda$ strictly increases.

It remains to prove that $R_\lambda$ has only three types of entries. Let
$\Omega=\binom{[N]}{2}$,
so that the rows and columns of $A_0$ are indexed by the two-element subsets of $[N]$, equivalently by the edges of $K_N$.
For $e,e'\in\Omega$, exactly one of the following three alternatives occurs:
\[
e=e',\quad |e\cap e'|=1,\quad |e\cap e'|=0.
\]
Let $A^{(0)},A^{(1)},A^{(2)}$ be the three $(0,1)$-matrices corresponding to these relations. Thus $A^{(0)}=I$, $A^{(1)}=A(L(K_N))=A_0$, and $A^{(2)}$ is the matrix whose $(e,e')$-entry is $1$ exactly when the two edges of $K_N$ are disjoint.
The span
\[
\mathcal B=\Span\{A^{(0)},A^{(1)},A^{(2)}\}
\]
is the Bose-Mesner algebra of the Johnson scheme $J(N,2)$.
Thus a matrix belongs to $\mathcal B$ if and only if its $(e,e')$-entry depends only on $|e\cap e'|$.

The algebra $\mathcal B$ is closed under multiplication.
Since $I$ and $A_0$ belong to $\mathcal B$, the matrix $\lambda I-A_0$ belongs to $\mathcal B$. Its inverse also belongs to $\mathcal B$: indeed, by the Cayley-Hamilton theorem, the inverse of an invertible matrix is a polynomial in that matrix, and $\mathcal B$ is closed under polynomials.
Consequently
\[
R_\lambda=(\lambda I-A_0)^{-1}\in\mathcal B.
\]
Therefore there exist numbers $\alpha,\beta,\gamma$ such that
\begin{equation}\label{abg}
 R_\lambda = \alpha A^{(0)} + \beta A^{(1)} + \gamma A^{(2)}= \alpha I + \beta A_0 + \gamma (J-I-A_0).
\end{equation}
Using the eigenspaces corresponding to \eqref{eq:triangular-spectrum}, and applying \eqref{abg} to an eigenvector of $A_0$ corresponding to the eigenvalue $N-4$ (respectively, the eigenvalue $-2$), we obtain
\begin{equation}\label{eq:alpha-beta-gamma}
\beta-\gamma=\frac1{(\lambda-N+4)(\lambda+2)}>0,
\quad
\alpha-\beta=\frac1{\lambda+2}+\beta-\gamma>0.
\end{equation}
Thus $\alpha>\beta>\gamma$.

For distinct triples $F,F'$, the value $b_F^T R_\lambda b_{F'}$ depends only on $|F\cap F'|$. Counting the nine ordered pairs of boundary edges gives
\[
b_F^\top R_\lambda b_{F'}=
\begin{cases}
\alpha+6\beta+2\gamma, & \text{if~} |F\cap F'|=2,\\
4\beta+5\gamma, & \text{if~} |F\cap F'|=1,\\
9\gamma, & \text{if~} |F\cap F'|=0.
\end{cases}
\]
After adding the $C_M$-term, the three possible off-diagonal entries of $C_M+B^\top R_\lambda B$ are
\[
w_2=1+\alpha+6\beta+2\gamma,
\quad
w_1=4\beta+5\gamma,
\quad
w_0=9\gamma.
\]
By \eqref{eq:alpha-beta-gamma},
\[
w_1-w_0=4(\beta-\gamma)>0,
\quad
w_2-w_1=1+(\alpha-\beta)+3(\beta-\gamma)>0.
\]
Hence $w_2>w_1>w_0$. The diagonal entries of $C_M+B^\top R_\lambda B$ are all equal to $3\alpha+6\beta$.

Let $M^*$ be a star family. For every $\lambda>\rho_0$, the matrix $C_M+B^T R_\lambda B$ is entrywise at most the corresponding matrix $C_{M^*}+B_{M^*}^T R_\lambda B_{M^*}$, because every off-diagonal entry in the latter has the largest possible weight $w_2$. Therefore
\[
\rho(C_M+B^\top R_\lambda B)\le
\rho(C_{M^*}+B_{M^*}^\top R_\lambda B_{M^*}).
\]
Evaluating at $\lambda^*=\rho(A_{M^*})$ and using \eqref{eq:schur-root} for $M^*$ gives
\[
\rho(C_M+B^\top R_{\lambda^*}B)\le \rho(C_{M^*}+B_{M^*}^\top R_{\lambda^*} B_{M^*}) = \lambda^*.
\]
Since the right side of \eqref{eq:schur-root} is strictly decreasing in $\lambda$ and the left side of \eqref{eq:schur-root} is strictly increasing, the solution for $M$ is at most $\lambda^*$. Thus $\q_1(\T_n^3\cup M)\le \q_1(\T_n^3\cup M^*)$.
\end{proof}

\begin{lemma}\label{lem:cospectral}
Let $M$ be a $3$-uniform family with $|M|=t\ge1$ such that any two distinct faces of $M$ share an edge. Then $\q_1(\T_n^3\cup M)$ depends only on $n$ and $t$, not on the isomorphism type of $M$.
\end{lemma}

\begin{proof}
The case $n=4$ is immediate, so assume $n\ge5$.
Use the same block matrix \eqref{eq:blockAM}. If $M$ is pairwise edge-intersecting, then $C_M=J_t-I_t$. Since the resolvent $R_\lambda$ belongs to the Bose-Mesner algebra of $J(N,2)$, the value $b_F^T R_\lambda b_{F'}$ depends only on $|F\cap F'|$. For distinct $F,F'\in M$, we always have $|F\cap F'|=2$. Hence for suitable functions $a(\lambda)$ and $b(\lambda)$ depending only on $n$,
\[
B^T R_\lambda B=a(\lambda)I_t+b(\lambda)(J_t-I_t).
\]
Thus the Schur complement
is
\[
\lambda I-C_M-B^TR_\lambda B= (\lambda-a(\lambda))I_t-(1+b(\lambda))(J_t-I_t).
\]
The matrix $C_M+B^TR_\lambda B$ is positive and has constant row sum
$a(\lambda)+(t-1)(1+b(\lambda))$, so its Perron vector is the all-ones vector.
By \eqref{eq:schur-root}, $\rho(A_M)$ is the unique solution $\lambda>\rho_0$ of
\begin{equation}\label{PerronRoot}
\lambda-a(\lambda)-(t-1)(1+b(\lambda))=0.
\end{equation}
This equation depends only on $n$ and $t$, and so does $\q_1(\T_n^3\cup M)=3+\rho(A_M)$.
\end{proof}

\begin{proof}[Proof of Theorem \ref{thm:main}]
Let $K\in \K(n,2,t)$ be extremal. By Theorem \ref{thm:structural}, $K=\T_n^3\cup M$ with $|M|=t$. Lemma \ref{lem:exact-comparison} shows that replacing $M$ by a star family $M^*$ cannot decrease the spectral radius. If $M$ is not pairwise edge-intersecting, put $H_M(\lambda)=C_M+B^TR_\lambda B$ and $\lambda^*=\rho(A_{M^*})$.
The comparison in Lemma \ref{lem:exact-comparison} gives $H_M(\lambda^*)\le H_{M^*}(\lambda^*)$ entrywise, with at least one strict inequality.
Both matrices are positive, so
$\rho(H_M(\lambda^*))<\rho(H_{M^*}(\lambda^*))=\lambda^*$.
The monotonicity in \eqref{eq:schur-root} then gives $\rho(A_M)<\lambda^*$, contradicting extremality.
Thus every extremal $M$ must be pairwise edge-intersecting.

Conversely, by Lemma \ref{lem:cospectral}, every pairwise edge-intersecting $M$ has the same spectral radius as the star family $M^*$ and is therefore extremal.

It remains to discuss uniqueness. For $t=1$, the statement is immediate. For $t=2$, two distinct triples must share an edge, and all such configurations are isomorphic. For $t\ge5$, use the following elementary fact: if a family of at least five triples is pairwise edge-intersecting, then all triples contain one common edge. Indeed, take two triples $A=\{a,b,c\}$ and $B=\{a,b,d\}$. Any further triple sharing an edge with both $A$ and $B$ either contains $\{a,b\}$ or is one of $\{a,c,d\}$ and $\{b,c,d\}$. If $\{a,c,d\}$ or $\{b,c,d\}$ belongs to the family, a triple $\{a,b,x\}$ sharing an edge with it must have $x\in\{c,d\}$. Thus the non-star case is contained in the four triples of $K_4^{(3)}$, and has size at most four. Hence for $t\ge5$ every pairwise edge-intersecting family is a star. For $t=3,4$, non-star subfamilies of $K_4^{(3)}$ are pairwise edge-intersecting and yield additional extremal complexes.
Indeed, for sufficiently large $n$, the apex lies in $\binom{n-1}{2}$ facets, while every other vertex lies in at most $n-2+t$ facets.
An isomorphism must therefore preserve the apex and the isomorphism type of $M$, so these additional complexes are not isomorphic to the star construction.
\end{proof}

\begin{proof}[Proof of Corollary \ref{cor:asymptotic}]
This is the asymptotic part of Theorem \ref{thm:structural}, together with Theorem \ref{thm:main}.
\end{proof}

\section{An extension to growing Betti numbers}

The Schur-complement comparison in Section \ref{sec:schur} is exact for the reduced class $\T_n^3\cup M$. The restriction to fixed $t$ comes from the structural reduction. Tracking the constants in Section \ref{sec:structure} gives the following range.

\begin{lemma}\label{lem:uniform-structural}
Let $t=t(n)\ge1$ be an integer satisfying $t=o(n^{1/4})$. If $K\in \K(n,2,t)$ is extremal, then for all sufficiently large $n$,
\[
K=\T_n^3\cup M,
\qquad
M\subseteq \binom{V\setminus\{u\}}3,
\qquad
|M|=t.
\]
\end{lemma}

\begin{proof}
The estimates in Lemmas \ref{lem:aligned} and \ref{lem:bounded-Mu} are uniform as long as the contradictions remain valid. Lemma \ref{lem:aligned} requires $n\gg t^2$ and gives $A_{\le1}=\emptyset$ and $A_2=A_2^u$. Lemma \ref{lem:bounded-Mu} then gives
\[
|M_u|<5t^2+10t.
\]
Consequently, for $F\in M_u$, \eqref{eq:small-noncone} gives
\[
f(\partial F)=O(t^2/n),
\]
and for every $x\in A$, \eqref{eq:ux-large} gives
\[
f(\{u,x\})>\frac12-O(t^4/n).
\]
If $t=o(n^{1/4})$, then $t^4/n=o(1)$ and $t^2/n=o(1)$. The face-switching argument in Lemma \ref{lem:cone-completion} therefore still gives a strict Rayleigh quotient increase whenever a cone face is missing.
Hence all cone faces are present, and the two bounds on $|M|$ in the proof of Theorem \ref{thm:structural} give $|M|=t$.
\end{proof}

\begin{proof}[Proof of Corollary \ref{cor:growing}]
By Lemma \ref{lem:uniform-structural}, every extremal complex has the form $\T_n^3\cup M$. Applying Lemma \ref{lem:exact-comparison} gives that a star family is extremal. The uniqueness statement follows from the same elementary classification of pairwise edge-intersecting triple families used in the proof of Theorem \ref{thm:main}.
\end{proof}

\begin{rmk}
The exact Schur-complement comparison shows that, after the cone reduction, the full range $1\le t\le n-3$ is settled. Thus the remaining obstacle to a full theorem for growing $t$ is the structural reduction. Improving Lemma \ref{lem:bounded-Mu}, or avoiding the loss $|M_u|^2/n$ in the final face-switching step, would enlarge the range of Corollary \ref{cor:growing}.
\end{rmk}

\end{document}